\pdfoutput=1
\documentclass[12pt,leqno]{amsart}
\usepackage{amssymb,amsmath,amsthm,bm}
\usepackage{graphicx}
\usepackage{color}
\usepackage[bookmarksopen,bookmarksdepth=3,colorlinks,citecolor=red,pagebackref,hypertexnames=true]{hyperref}
\usepackage[msc-links,nobysame,sorted-cites,initials]{amsrefs}
\usepackage[inline]{enumitem}
\setlist[enumerate]{font=\upshape}
\usepackage{mathtools}
\mathtoolsset{showonlyrefs}
\definecolor{darkblue}{RGB}{0,0,160}

\usepackage[lining]{libertine}   
\usepackage{cabin}
\usepackage[libertine]{newtxmath}

\def\eps{\varepsilon} 
\def\d{{\mathrm d}}
\def\dist{\mathrm {dist}}

\def\R {\mathbb{R}}

\def\M {{\mathrm M}}
\def\Z {{\mathbb Z}}

\def\P {\mathbb P}

\def\ind{\mathbf{1}}

\newcommand{\supp}{\operatorname{supp}}

\newcommand{\PM}{\mathsf{P}}

\newcommand{\nat}[1]{\langle #1\rangle^{\natural}}
\newcommand{\natl}[2]{\langle #1\rangle^{\natural,[#2]}}

   \def\XXint#1#2#3{{\setbox0=\hbox{$#1{#2#3}{\int}$}
       \vcenter{\hbox{$#2#3$}}\kern-.5\wd0}}

\makeatletter
\DeclareRobustCommand\widecheck[1]{{\mathpalette\@widecheck{#1}}}
\def\@widecheck#1#2{%
    \setbox\z@\hbox{\m@th$#1#2$}%
    \setbox\tw@\hbox{\m@th$#1%
       \widehat{%
          \vrule\@width\z@\@height\ht\z@
          \vrule\@height\z@\@width\wd\z@}$}%
    \dp\tw@-\ht\z@
    \@tempdima\ht\z@ \advance\@tempdima2\ht\tw@ \divide\@tempdima\thr@@
    \setbox\tw@\hbox{%
       \raise\@tempdima\hbox{\scalebox{1}[-1]{\lower\@tempdima\box
\tw@}}}%
    {\ooalign{\box\tw@ \cr \box\z@}}}
\makeatother

\newcounter{counter}

\numberwithin{equation}{section}
\numberwithin{counter2}{section}

\newtheorem{theorem}[counter]{Theorem}

\numberwithin{withinsub}{subsection}
\newtheorem{corollary}{Corollary}
\numberwithin{corollary}{counter}
\newtheorem{lemma}[subsection]{Lemma} 
 
\newtheorem{proposition}[subsection]{Proposition}

\theoremstyle{definition}

\newtheorem*{remark*}{Remark}
\newtheorem*{warn*}{A word of warning}

\newtheorem{remark}[subsection]{Remark} 
 
\theoremstyle{definition}

 \theoremstyle{plain}

\numberwithin{figure}{section}

\begin{document}

\title[Percentile oscillations and the intrinsic square function]{Percentile oscillations and the intrinsic square function}

  \author[F.\ Di Plinio]{Francesco Di Plinio}
\address[F.\ Di Plinio]{Dipartimento di Matematica e Applicazioni ``R.\ Caccioppoli'' \\ \newline\indent Universit\`a di Napoli ``Federico II'', Via Cintia, Monte S.\ Angelo 80126 Napoli, Italy}
\email{\href{mailto:francesco.diplinio@unina.it}{\textnormal{francesco.diplinio@unina.it}}}


\begin{abstract}
Wilson's intrinsic square function tests the input against all the cancellative functions of a given smoothness, at every location and scale. We show that it is dominated pointwise by a   sparse square function, whose coefficients are the percentile oscillations of the input modulo polynomials, at every polynomial order below the smoothness plus one. These coefficients are local sizes, which nevertheless retain essentially all the cancellation: the weighted norms of the sparse square function are equivalent to smooth grand maximal norms, for all $0<p<\infty$ and $A_\infty$ weights, with explicit dependence on the weight. These are weighted Hardy norms in the range allowed by the smoothness. Consequences include potential bounds, and weighted Hardy space estimates with the sharp factor $[w]_{A_\infty}^{1/2}$. For the derivatives of Calder\'on--Zygmund operators, we characterize the corresponding square function estimates by testing conditions on the polynomial symbols of their wavelet resolution. In the quadratic setting, the factor $[w]_{A_\infty}^{1/2}$ detects a transition from $\mathrm{BMO}$ testing conditions to vanishing symbols.
\end{abstract}
\maketitle

\section{Introduction}\label{s:intro}
Square functions measure the cancellation of a function across scales. Among them, the intrinsic square function of Wilson \cite{Wilson07} is universal, in the sense that at each location and scale, it tests the input against all the cancellative functions of a given smoothness at once. Thus it dominates the classical area functions with smooth cancellative kernels. It also characterizes the weighted Hardy spaces of \cite{ST89}, in a range of exponents and weights depending on the smoothness \cite{WL12}. 
This article shows that, at every polynomial order below the smoothness plus one, the intrinsic square function admits a sparse version, which retains essentially all of its cancellation. Its coefficients are the \emph{percentile oscillations} of the input, introduced in \cite{DPModPoly} as an extension of the coefficients of order zero of Conde-Alonso, Lorist and Rey \cite{CLR26}. 

The percentile oscillation $\nat{f}_{m,Q}$ of order $m\geq0$ of a function $f$ on a cube $Q$ is a size, rather than a pairing with a cancellative function. It is a fixed upper percentile on $Q$ of the localized smooth grand maximal function of $f-\pi_Q$, where $\pi_Q$ is the polynomial approximation of degree less than $m$ of $f$ on $3Q$, see \eqref{e:po}. Thus it depends only on $f$ on $3Q$, it vanishes exactly when $f$ coincides there with a polynomial of degree less than $m$, and it sees the cancellation of $f$ at all the scales below that of $Q$. 

 Theorem \ref{t:square}, the first main result of this work, shows that the intrinsic square function is pointwise dominated  by a quadratic sparse operator whose coefficients are the percentile oscillations,  at each order below the smoothness plus one.  In turn, weighted $L^p$ norms of the quadratic sparse operator are equivalent to smooth grand maximal norms, for all $0<p<\infty$ and all $A_\infty$ weights, with explicit dependence on the weight, see \eqref{e:sparsechar}. Their identification with weighted Hardy norms requires the usual range restrictions on $p$, $w$ and the smoothness.

The order $m$ adapts the percentile oscillations to Sobolev spaces, as $\nat{f}_{m,Q}$ is at most a constant multiple of $\ell(Q)^m$ times the percentile oscillation of order zero of $\nabla^mf$ \cite{DPModPoly}*{Proposition 6.5}. Thus Theorem \ref{t:square} leads to potential bounds, Corollary \ref{c:sqpot}, and to weighted bounds in terms of the maximal percentile oscillation $\mathcal N_mf$, with the sharp factor $[w]_{A_\infty}^{1/2}$, for all $0<p<\infty$, Corollary \ref{c:weighted}. Theorem \ref{t:square} recovers, through \eqref{e:poelem}, the sharp weighted estimates for intrinsic square functions of Lerner \cite{lerner11}, and, for $m=0$, it is the Euclidean counterpart of the cancellative quadratic sparse bounds of \cite{CLR26}. 

The order $m$ also  appears in the dependence on the weight, which distinguishes the testing conditions of the weighted Hardy--Sobolev theory of the derivatives of Calder\'on--Zygmund operators \cite{DPModPoly}. Theorem \ref{t:SoT} characterizes instead   square function estimates for the derivatives $\partial^\alpha T$, $|\alpha|=\kappa$, of a Calder\'on--Zygmund operator $T$, in terms of the maximal percentile oscillations of the derivatives of order $\kappa$ of the input, by testing conditions on the polynomial symbols of the wavelet resolution of $T$ \cite{DGWdom}. With $\mathcal N_0(\partial^\beta f)$ on the right, that is, at residual order zero, the lower symbols must vanish, the top symbols must belong to $\mathrm{BMO}$, and $T^*1$ must be constant when $\kappa=0$, while with $\mathcal N_1(\partial^\beta f)$, at residual order one, all the symbols must vanish. Thus, as the estimates become more demanding, the conditions on the top symbols pass through
$
\mathrm{BMO}\longrightarrow L^\infty\longrightarrow 0,
$
corresponding respectively to residual order zero, positive residual order with linear $A_\infty$ dependence \cite{DPModPoly}*{Theorem A}, and residual order one in the quadratic theory, with the sharp factor $[w]_{A_\infty}^{1/2}$.

The proof of Theorem \ref{t:square} combines the stopping principle of \cite{DPModPoly}*{Proposition 5.1} with a balayage argument, Proposition \ref{c:balayage}. As the intrinsic coefficients of a cube see the whole input, and not only its restriction to a neighborhood of the cube, the local bounds have tails, which are charged to the stopping cubes and to the dilates of the top cube containing a given point. Theorem \ref{t:SoT} is proved through the wavelet representation of \cite{DPModPoly}*{Proposition 7.6}, and the vanishing of the top symbols in its part (b) by testing on the logarithmic functions and the power weights of \cite{DPModPoly}*{\S11.3}. The main definitions of \cite{DPModPoly} are recalled in \S\ref{ss:notation} and \S\ref{ss:SoTmain}, and its numbered statements are cited with precise references.

\section{Main results}\label{s:main}
\subsection{Notation}\label{ss:notation}
The notation is that of \cite{DPModPoly}*{Section 2}.  We limit ourselves to discussing the  points which are nonstandard or not completely intuitive. Generic constants $C$ and implicit constants  in $\lesssim $ are termed \emph{admissible} when they depend only on $d$, on the polynomial order $m$ and on the parameters of the classes below. The sign $\lesssim_a$ indicates that the implicit constant also depends on the parameter $a$. Hereafter, $\M$ is the Hardy--Littlewood maximal function; $\M_s$ is its fractional version of order $s$. We use the Fujii--Wilson characteristic $[w]_{A_\infty}\coloneqq\sup_Qw(Q)^{-1}\int_Q\M(w\ind_Q)$. A collection $\mathcal S$ is $\eta$-sparse if it admits pairwise disjoint measurable sets $G_S\subseteq S$ with $|G_S|\geq\eta|S|$.

For $\sigma>0$, $\rho\geq0$ and a cube $R$, the \emph{wavelet class} $\Phi^{\sigma,\rho}(R)$ consists of the $L^1$-normalized rescalings
\[
\varphi_R\coloneqq\frac1{|R|}\,\varphi\left(\frac{\cdot-c(R)}{\ell(R)}\right)
\]
of the functions $\varphi$ such that, for $|\alpha|\leq\lfloor\sigma\rfloor$, $x\in\R^d$ and $0<|h|\leq1$,
\begin{equation}\label{e:classes}
(1+|x|)^{d+\sigma+\rho}\left(\left|\partial^\alpha\varphi(x)\right|+\frac{\left|\partial^\alpha\varphi(x+h)-\partial^\alpha\varphi(x)\right|}{|h|^{\sigma-\lfloor\sigma\rfloor}}\right)\leq1,
\end{equation}
and the \emph{cancellative class} $\Psi^{\sigma,\rho}(R)$ consists of the elements of $\Phi^{\sigma,\rho}(R)$ with vanishing moments of order at most $\lfloor\sigma\rfloor$. We write \[\varphi(f)\coloneqq\int\varphi f,\] and $\varphi\in C\Phi^{\sigma,\rho}(R)$ if $\varphi$ is an admissible multiple of an element of $\Phi^{\sigma,\rho}(R)$; likewise for the other classes.

Fix $\tau>0$. The localized smooth grand maximal function of \cite{CLR26} is
\begin{equation}\label{e:Mtau}
\M^\tau_Qf(x)\coloneqq\sup\left\{|\varphi(f)|:\;\varphi\in\Phi^{\tau,1}(R),\;\supp\varphi\subseteq R,\;x\in R\subseteq Q\right\},\qquad x\in Q,
\end{equation}
and $\M^\tau_{\R^d}$ is defined with $\R^d$ in place of $Q$. For $m\geq1$, $\P^{m-1}$ is the space of the polynomials of degree less than $m$, and the averaged Taylor polynomial of \cite{DPModPoly}*{(2.5)},
\[
\PM^m_Qf(x)\coloneqq\sum_{|\alpha|<m}\frac1{\alpha!}\int_{\R^d}\theta_Q(y)\,\partial^\alpha f(y)\,(x-y)^\alpha\,\d y,
\]
where $\theta_Q$ is the $L^1$-normalized rescaling of a fixed smooth nonnegative bump $\theta$ supported in $(-\frac18,\frac18)^d$ with $\int\theta=1$, extends by integration by parts to a projection of $L^1_{\mathrm{loc}}(\R^d)$ onto $\P^{m-1}$; also, $\P^{-1}\coloneqq\{0\}$ and $\PM^0_Q\coloneqq0$. With the rearrangement notation \[
g^{*,Q}(r)\coloneqq   \inf \left\{s>0:|\{x\in Q:|g(x)|>s\}|\leq r|Q|\right\},
\] the \emph{percentile oscillation} of order $m\geq0$ of $f\in L^1_{\mathrm{loc}}(\R^d)$ on the cube $Q$ is
\begin{equation}\label{e:po}
\nat{f}_{m,Q}\coloneqq\left(\M^\tau_{3Q}\big[f-\PM^m_{3Q}f\big]\right)^{*,Q}(r),
\end{equation}
which for $m=0$ is the coefficient of \cite{CLR26}, and the \emph{maximal percentile oscillation} is \[\mathcal N_mf\coloneqq\sup_Q\nat{f}_{m,Q}\ind_Q,\] over all cubes $Q$. As in \cite{DPModPoly}*{\S2.3}, $\tau$ and $r$ are omitted from the notation, $0<r\leq r_0$ for a small admissible $r_0$, which may decrease from statement to statement, implied constants do not depend on $r$ unless indicated, and $\natl{f}{r}_{m,Q}$, $\mathcal N^{[r]}_m$ indicate the level. By \cite{DPModPoly}*{Lemma 2.4(ii)--(iii)},
\begin{equation}\label{e:poelem}
\nat{f}_{m,Q}\lesssim\inf_{\pi\in\P^{m-1}}\nat{f-\pi}_{0,Q},\qquad \nat{f}_{m,Q}\lesssim\frac1r\inf_{\pi\in\P^{m-1}}\left\langle|f-\pi|\right\rangle_{3Q},
\end{equation}
so that $\mathcal N_nf\lesssim\mathcal N_mf$ for $n\geq m$, and $\mathcal N_mf\lesssim r^{-1}\M f$. As customary, $\mathrm I_s$ is the Riesz potential of order $s$,
\[
\mathrm I_sg(x)\coloneqq\int_{\R^d}\frac{g(y)}{|x-y|^{d-s}}\,\d y, \qquad x \in \mathbb R^d.
\]

\subsection{The intrinsic square function}\label{ss:sqmain}
The \emph{intrinsic square function}, a version of that of Wilson \cite{Wilson07} built on the cancellative classes, is
\begin{equation}\label{e:intrinsic}
\mathrm S^\tau_{\mathcal R}f\coloneqq\Big(\sum_{R\in\mathcal R}\left|\left\langle f,\Psi^{\tau,1}(R)\right\rangle\right|^2\ind_R\Big)^{\frac12},
\end{equation}
defined for $f\in L^2(\R^d)$ and any subcollection $\mathcal R$ of a dyadic grid, where, for a class $\Phi$ of test functions,
\[
\left|\left\langle f,\Phi\right\rangle\right|\coloneqq\sup_{\varphi\in\Phi}|\varphi(f)|.
\] We write $\mathrm S^\tau\coloneqq\mathrm S^\tau_{\mathcal D}$; other dyadic grids lead to equivalent square functions, see \eqref{e:gridchange}. Below, $\mathrm S^\tau$, $\nat{f}_{m,S}$ and $\mathcal N_m$ are computed with the same smoothness parameter $\tau$. The first main result is a quadratic sparse estimate for the square function  $\mathrm S^\tau$.
\begin{theorem}\label{t:square}
Let $m\geq0$ be an integer with $\tau>m-1$, and let $0<r\leq r_0$. For each $f\in L^2(\R^d)$ with compact support there exists a sparse collection $\mathcal S$ of cubes of $\R^d$ such that, a.e.\ on $\R^d$,
\[
\mathrm S^\tau f\lesssim\Big(\sum_{S\in\mathcal S}\left[\nat{f}_{m,S}\right]^2\ind_S\Big)^{\frac12},\qquad\M^\tau_{\R^d}f\lesssim\sum_{S\in\mathcal S}\nat{f}_{m,S}\ind_S.
\]
\end{theorem}
The second estimate is the linear sparse bound of \cite{DPModPoly}*{Proposition 9.5}, here with the same collection as the first.
\begin{corollary}\label{c:sqpot}
Let $m\geq1$ be an integer with $\tau>m-1$, and let $0<r\leq r_0$. For every $f\in \mathcal{C}^\infty_c(\R^d)$ there exists a sparse collection $\mathcal S$ of cubes of $\R^d$ such that, a.e.\ on $\R^d$,
\begin{equation}\label{e:sqfrac}
\mathrm S^\tau f\lesssim\Big(\sum_{S\in\mathcal S}\Big[\ell(S)^m\left(\M^\tau_{3S}\left[\nabla^mf\right]\right)^{*,S}(r)\Big]^2\ind_S\Big)^{\frac12}.
\end{equation}
Consequently, a.e.\ on $\R^d$,
\[
\mathrm S^\tau f\lesssim\frac1r\,\mathrm I_m\left(|\nabla^mf|\right),\quad m<d,\qquad \mathrm S^\tau f\lesssim\frac1r\left[\mathrm I_{2m}\left(|\nabla^mf|^2\right)\right]^{\frac12},\quad 2m<d.
\]
\end{corollary}
There is no general pointwise ordering between the right hand sides of the two potential bounds, and both discard the cancellation retained by \eqref{e:sqfrac}.
\begin{corollary}\label{c:weighted}
Let $m\geq0$ be an integer with $\tau>m-1$, and let $0<r\leq r_0$, $0<p<\infty$ and $w\in A_\infty$. For every $f\in L^2(\R^d)$ with $\mathcal N_mf\in L^p(w)$,
\begin{equation}\label{e:weighted}
\left\|\mathrm S^\tau f\right\|_{L^p(w)}\lesssim_p[w]_{A_\infty}^{\frac12}\left\|\mathcal N_mf\right\|_{L^p(w)}.
\end{equation}
The implied constant grows at most like $\sqrt p$ as $p\to\infty$. If $1\leq m<d$, then, for every $f\in \mathcal{C}^\infty_c(\R^d)$ with $\M_m(|\nabla^mf|)\in L^p(w)$, \eqref{e:weighted} holds with $\mathcal N_mf$ replaced by $r^{-1}\M_m(|\nabla^mf|)$.
\end{corollary}
The factor $[w]_{A_\infty}^{1/2}$ and the growth $\sqrt p$ in \eqref{e:weighted} cannot be improved, see Remark \ref{r:sharp}. Together, Corollary \ref{c:weighted} and \cite{DPModPoly}*{Proposition 2.8} control the square and grand maximal norms by the maximal percentile oscillation, and hence the corresponding weighted Hardy norms in their admissible range, with the sharp factors $[w]_{A_\infty}^{1/2}$ for the square function norm, and $[w]_{A_\infty}$ for the maximal norm, the latter for $m\geq1$ \cite{DPModPoly}*{Remark 2.9}.

The quadratic sparse operators of Theorem \ref{t:square} characterize the smooth grand maximal norms. For fixed $p$ and $w$, their norms in $L^p(w)$ are comparable to those of $\mathcal N_mf$ and of $\M^\tau_{\R^d}f$, with the value $+\infty$ allowed. Indeed, let $m$, $p$ and $w$ be as in Corollary \ref{c:weighted}, so that $w\in A_q$ for some $1<q<\infty$, and let $f\in L^2(\R^d)$ have compact support. Then the collection $\mathcal S$ of Theorem \ref{t:square} satisfies
\begin{equation}\label{e:sparsechar}
\begin{split}
\big\|\M^\tau_{\R^d}f\big\|_{L^p(w)}&\lesssim\left\|\sum_{S\in\mathcal S}\nat{f}_{m,S}\ind_S\right\|_{L^p(w)}\\&\lesssim_p[w]_{A_\infty}\left\|\left(\sum_{S\in\mathcal S}\left[\nat{f}_{m,S}\right]^2\ind_S\right)^{\frac12}\right\|_{L^p(w)}\\&\lesssim_p[w]_{A_\infty}^{\frac32}\left\|\mathcal N_mf\right\|_{L^p(w)}\\&\lesssim[w]_{A_\infty}^{\frac32}\left(\frac Cr\right)^{\frac qp}[w]_{A_q}^{\frac1p}\big\|\M^\tau_{\R^d}f\big\|_{L^p(w)},
\end{split}
\end{equation}
where $C$ depends only on $d$. The first inequality is the second estimate of Theorem \ref{t:square}, while the second follows from \cite{DPModPoly}*{Lemma 9.3}, as the supremum of a sequence is at most its $\ell^2$ norm. The third bound originates from Lemma \ref{l:sqAinf}, as $\nat{f}_{m,S}\ind_S\leq\mathcal N_mf$, while the fourth is Lemma \ref{l:median}. Thus the coefficients $\nat{f}_{m,S}$, which vanish when $f$ coincides on $3S$ with a polynomial of degree less than $m$, still control $\M^\tau_{\R^d}f$, also for $p\leq1$, where the sharp maximal function of a nonzero compactly supported function does not belong to $L^p(\R^d)$. Moreover, $\mathrm S^\tau f$ dominates pointwise, up to a constant, the Lusin area function of $f$ built on any Schwartz function with vanishing moments up to order $\lfloor\tau\rfloor$. As those with nondegenerate kernels characterize the weighted Hardy spaces $H^p_w$ of \cite{ST89} in a range of $p$ and $w$ depending on $\tau$, so does $\mathrm S^\tau$ where $\M^\tau_{\R^d}$ does, with the sharp factor $[w]_{A_\infty}^{1/2}$ in \eqref{e:weighted}.

\subsection{Square functions of Calder\'on--Zygmund operators}\label{ss:SoTmain}
Let $k\geq0$ be an integer and $0<\delta\leq1$. Informally, a \emph{Calder\'on--Zygmund operator of order $k$} is an operator $T$ with the weak boundedness property, acting on spaces of test functions defined through their wavelet coefficients, whose kernel $K$ satisfies, for $0\leq i\leq k$,
\[
\begin{split}
\left|\nabla^i_xK(x,y)\right|&\leq\frac1{|x-y|^{d+i}},\\
\left|\nabla^i_xK(x',y')-\nabla^i_xK(x,y)\right|&\leq\frac{(|x-x'|+|y-y'|)^\delta}{|x-y|^{d+i+\delta}},\quad |x-x'|+|y-y'|\leq\tfrac12|x-y|;
\end{split}
\]
see \cite{DPModPoly}*{Definition 3.3} for the precise definition. The kernel and weak boundedness constants are normalized to one. Its derivatives $\partial^\alpha T$, $|\alpha|=\kappa\leq k$, are defined by duality,
\[
\langle\partial^\alpha Tf,g\rangle\coloneqq(-1)^\kappa\langle Tf,\partial^\alpha g\rangle.
\]

We use the compactly supported wavelets of Daubechies \cite{Dau88}, as in \cite{DGWdom}*{Proposition 2.3 and Remark 2.2} and in \cite{DPModPoly}*{\S3.1}. There are real valued functions $\psi_S$, $2^d-1$ for each $S\in\mathcal D$, such that
\[
\psi_S\in C\Psi^{k+1,1}(S),\qquad\supp\psi_S\subseteq CS,
\]
and $\{|S|^{\frac12}\psi_S:S\in\mathcal D\}$ is an orthonormal basis of $L^2(\R^d)$. Thus, with sums over $S$ including the index of the wavelet, and with the pairing $\langle\cdot,\cdot\rangle$ bilinear,
\begin{equation}\label{e:onb}
f=\sum_{S\in\mathcal D}|S|\,\psi_S(f)\,\psi_S,\qquad\langle f,g\rangle=\sum_{S\in\mathcal D}|S|\,\psi_S(f)\,\psi_S(g),\qquad f,g\in L^2(\R^d).
\end{equation}
The \emph{symbols} of order $\kappa$ of $T$ are the sequences
\[
b^{\alpha,\gamma}_S\coloneqq\left\langle\partial^\alpha T\left[(\cdot-c(S))^\gamma\right],\psi_S\right\rangle,\quad |\gamma|\leq\kappa=|\alpha|,\qquad b^{0,\star}_S\coloneqq\langle T\psi_S,\mathbf 1\rangle,\qquad S\in\mathcal D.
\]
See \cite{DPModPoly}*{(3.5), (3.6) and Definition 3.4} for the rigorous definitions. Thus $b^{0,0}$ and $b^{0,\star}$ are the wavelet coefficients of $T1$ and $T^*1$. The symbols $b^{\alpha,\gamma}$ with $|\gamma|<\kappa$ and with $|\gamma|=\kappa$ are the \emph{lower} and the \emph{top} symbols, respectively. A symbol $b$ \emph{vanishes} if $b_S=0$ for all $S$, and \emph{belongs to} $\mathrm{BMO}(\R^d)$, respectively $L^\infty(\R^d)$, if $b_S=\psi_S(\mathfrak b)$ for some $\mathfrak b$ in $\mathrm{BMO}(\R^d)$, respectively $L^\infty(\R^d)$, unique modulo constants; in the first case, the norm $\|b\|_{\mathrm{BMO}}$ is that of $\mathfrak b$. Theorem \ref{t:SoT}, the quadratic counterpart of \cite{DPModPoly}*{Theorem A}, characterizes the operators for which the square functions of $\partial^\alpha Tf$ are controlled by the maximal percentile oscillations of the derivatives of order $\kappa$ of $f$.
\begin{theorem}\label{t:SoT}
Let $T$ be a Calder\'on--Zygmund operator of order $k$, $0\leq\kappa\leq k$, $|\alpha|=\kappa$, $0<\tau<\delta$ and $0<r\leq r_0$. Consider the square function estimates
\begin{align}
\label{e:SoT}
\left\|\mathrm S^\tau\left(\partial^\alpha Tf\right)\right\|_{L^p(w)}&\lesssim_p[w]_{A_\infty}^{\frac12}\sum_{|\beta|=\kappa}\left\|\mathcal N_0\left(\partial^\beta f\right)\right\|_{L^p(w)},
\\ \label{e:SoT1}
\left\|\mathrm S^\tau\left(\partial^\alpha Tf\right)\right\|_{L^p(w)}&\lesssim_p[w]_{A_\infty}^{\frac12}\sum_{|\beta|=\kappa}\left\|\mathcal N_1\left(\partial^\beta f\right)\right\|_{L^p(w)},
\end{align}
each required for all $f\in\mathcal{C}^\infty_c(\R^d)$, $0<p<\infty$ and $w\in A_\infty$; part of each requirement is that $\partial^\alpha Tf\in L^2(\R^d)$.
\begin{enumerate}[label=(\alph*)]
\item The estimate \eqref{e:SoT} holds if and only if the lower symbols vanish, the top symbols belong to $\mathrm{BMO}(\R^d)$, and, if $\kappa=0$, the symbol $b^{0,\star}$ vanishes. Its implied constant is then ${\lesssim_p1+\max_{|\gamma|=\kappa}\|b^{\alpha,\gamma}\|_{\mathrm{BMO}}}$.
\item The estimate \eqref{e:SoT1} holds if and only if the lower and the top symbols vanish, and, if $\kappa=0$, so does $b^{0,\star}$. Its implied constant is then $\lesssim_p1$.
\end{enumerate}
In both cases, the implied constants grow at most like $\sqrt p$ as $p\to\infty$, and the factor $[w]_{A_\infty}^{1/2}$ cannot be improved in general.
\end{theorem}
The \emph{residual order} of the estimates \eqref{e:SoT} and \eqref{e:SoT1}, and of the weighted Hardy--Sobolev inequalities \cite{DPModPoly}*{(3.9)}, which \cite{DPModPoly}*{Theorem A} characterizes with the factor $[w]_{A_\infty}$, is the order of the maximal percentile oscillations on their right hand sides: zero, one and $m-\kappa$, respectively. At residual order zero, for $\kappa\geq1$, the conditions in (a) coincide with the conditions (i) of \cite{DPModPoly}*{Theorem A} for $m=\kappa$. For $\kappa=0$, the conditions in (a) say that $T$ is bounded on $L^2(\R^d)$, by the $T1$ theorem \cite{DJ84}, and that $T^*1=0$ in $\mathrm{BMO}$, whereas \cite{DPModPoly}*{Theorem A} only requires $T^*1\in\mathrm{BMO}(\R^d)$. The vanishing of $T^*1$ is necessary, as the square function of an integrable function with nonvanishing mean is not integrable at infinity.

At residual order one, the conditions in (b) require the top symbols to vanish, whereas the conditions (i) of \cite{DPModPoly}*{Theorem A} for $m=\kappa+1$ only require them to belong to $L^\infty(\R^d)$. Top symbols in $L^\infty(\R^d)$ do not suffice for (b). Indeed, if $b\in L^\infty(\R^d)$ is nonconstant, the canonical paraproduct $P_b$ of \cite{DPModPoly}*{(7.8)}, whose averaging functions $\vartheta_Q$ are supported in $Q$, with $\int\vartheta_Q=1$, is up to a constant multiple a Calder\'on--Zygmund operator of order zero with $b^{0,0}_S=\psi_S(b)$ and $b^{0,\star}=0$ \cite{DPModPoly}*{Lemma 10.7(iii)}. Thus $P_b$ satisfies \cite{DPModPoly}*{(3.9)} for every $m\geq0$, by \cite{DPModPoly}*{Theorem A}, but not \eqref{e:SoT1}, by (b), as $\psi_S(b)\neq0$ for some $S$. Such a paraproduct reads the averages of the top derivatives of $f$, which may carry a full logarithm along the scales: the linear dependence on $[w]_{A_\infty}$ absorbs it, while the factor $[w]_{A_\infty}^{1/2}$ only absorbs its square root.

\section{A balayage principle}\label{s:balayage}
The next proposition is a balayage version of the stopping principle of \cite{DPModPoly}*{Proposition 5.1}, in which the local weak-type assumption is allowed to have tails, measured by the local bounds of the larger cubes whose good sets are close to the given cube. These cubes range in the collection
\[
\mathcal D^+(Q_0)\coloneqq\mathcal D(Q_0)\cup\left\{2^jQ_0:j\geq1\right\}.
\]
Precisely, with $\eps_d\coloneqq2^{-d-4}$ as in \cite{DPModPoly}*{Proposition 5.1}, given measurable sets $E_P\subseteq P$ with $|E_P|\leq\eps_d|P|$ and numbers $\lambda_P\geq0$, $P\in\mathcal D^+(Q_0)$, define the \emph{tailed local bounds}
\begin{align}
\label{e:tailed}
\lambda^\dagger_Q&\coloneqq\lambda_Q+\sum_{k\geq1}2^{-k}\inf\left\{\lambda_P:P\in\mathcal T_k(Q)\right\},\qquad Q\in\mathcal D(Q_0),
\\ \nonumber
\mathcal T_k(Q)&\coloneqq\left\{P\in\mathcal D^+(Q_0):P\supseteq Q,\,\ell(P)\geq2^k\ell(Q),\,3\cdot2^kQ\cap(P\setminus E_P)\neq\varnothing\right\}.
\end{align}
Each $\mathcal T_k(Q)$ is nonempty, as it contains the dyadic ancestor of $Q$ of sidelength $2^k\ell(Q)$ when $2^k\ell(Q)\leq\ell(Q_0)$, and the smallest cube $2^jQ_0$ of sidelength at least $2^k\ell(Q)$ otherwise. Without tails, \eqref{e:P4} below is a pointwise form of the packing condition of the Carleson sequences of \cite{DPZyg}, and \eqref{e:tailed} is the analogue of the tailed averages of \cite{DPRDF}*{\S2.2}.
\begin{proposition}\label{c:balayage}
Let $Q_0$ be a cube and $K\geq1$. For $P\in\mathcal D^+(Q_0)$, let $E_P\subseteq P$ be measurable with $|E_P|\leq\eps_d|P|$, and $\lambda_P\geq0$, such that the tailed local bounds \eqref{e:tailed} are finite. Let $\{a_R:R\in\mathcal D(Q_0)\}$ be nonnegative numbers such that
\begin{equation}\label{e:P4}
\sum_{\substack{R\in\mathcal D(Q)\\ R\ni y}}a_R\leq K\lambda^\dagger_Q,\qquad Q\in\mathcal D(Q_0),\; y\in Q\setminus E_Q.
\end{equation}
Then there exists a $\frac12$-sparse collection $\mathcal S\subseteq\mathcal D^+(Q_0)$, containing the cubes $2^jQ_0$, $j\geq1$, such that
\begin{equation}\label{e:principle2}
\sum_{R\in\mathcal D(Q_0)}a_R\ind_R\lesssim K\sum_{S\in\mathcal S}\lambda_S\ind_S\qquad\textrm{a.e.\ on } \R^d.
\end{equation}
\end{proposition}
\begin{proof}
Apply \cite{DPModPoly}*{Proposition 5.1} to
\[
F_Q^{(n)}\coloneqq\sum_{\substack{R\in\mathcal D(Q)\\\ell(R)\geq2^{-n}}}a_R\ind_R,
\qquad Q\in\mathcal D(Q_0),
\]
with exceptional sets $E_Q$ and bounds $\lambda_Q^\dagger$. The local bound is \eqref{e:P4}; the oscillation condition holds because $F_Q^{(n)}-F_L^{(n)}$ is constant on every proper dyadic subcube $L$ of $Q$. The resulting $\frac34$-sparse collection $\mathcal S_0$ contains $Q_0$, depends only on the sets $E_Q$, and satisfies \cite{DPModPoly}*{(5.5)}. Thus, letting $n\to\infty$,
\begin{equation}\label{e:balstep}
\sum_{R\in\mathcal D(Q_0)}a_R\ind_R\leq3K\sum_{S\in\mathcal S_0}\lambda_S^\dagger\ind_S
\qquad\textrm{a.e. on }Q_0.
\end{equation}
Set $\mathcal S\coloneqq\mathcal S_0\cup\{2^jQ_0:j\geq1\}$. The annuli $2^jQ_0\setminus2^{j-1}Q_0$ witness its $\frac12$-sparseness outside $Q_0$.

Fix $x\in Q_0$ belonging to finitely many cubes of $\mathcal S_0$; this holds almost everywhere, since the $n$-th stopping generation has measure at most $4^{-n}|Q_0|$. Write $\mathcal C_x\coloneqq\{P\in\mathcal S:x\in P\}$. For $S\in\mathcal C_x\cap\mathcal S_0$ and $k\geq1$, let $P=P(S,k)$ be the least cube of $\mathcal C_x$ with $\ell(P)\geq2^k\ell(S)$, and let $P^-$ be its predecessor. Then
\begin{equation}\label{e:balchain}
S\subseteq P^-,\qquad \ell(P^-)<2^k\ell(S),\qquad P\in\mathcal T_k(S).
\end{equation}
For the last assertion, if $P\in\mathcal S_0$, then $3\cdot2^kS\supseteq P^-$ and $|P^-\cap E_P|\leq\frac14|P^-|$, by the second inequality in \cite{DPModPoly}*{(5.5)}, as $P^-$ is maximal among the elements of $\mathcal S_0$ strictly contained in $P$: any such element containing $P^-$ contains $x$. Otherwise $\ell(P)=2\ell(P^-)<2^{k+1}\ell(S)$, so that $3\cdot2^kS$ meets at least $2^{-d}|P|>|E_P|$ of $P$.

For fixed $S$ and $P$, the indices assigned to $P$ satisfy
\[
\sum_{k:\,P(S,k)=P}2^{-k}\leq2\frac{\ell(S)}{\ell(P^-)}.
\]
The sidelengths along $\mathcal C_x$ at least double. Consequently,
\[
\begin{split}
\sum_{S\in\mathcal C_x\cap\mathcal S_0}\lambda_S^\dagger
&\leq\sum_{S\in\mathcal C_x\cap\mathcal S_0}\lambda_S
+2\sum_{P\in\mathcal C_x}\lambda_P
  \sum_{\substack{S\in\mathcal C_x\cap\mathcal S_0\\S\subseteq P^-}}
  \frac{\ell(S)}{\ell(P^-)}
\leq5\sum_{P\in\mathcal C_x}\lambda_P.
\end{split}
\]
Together with \eqref{e:balstep}, this proves the claim on $Q_0$. Outside $Q_0$ the left hand side of \eqref{e:principle2} vanishes.
\end{proof}
\section{Sparse domination of the intrinsic square function}\label{s:square}
This section contains the proof of Theorem \ref{t:square}, the quadratic counterpart of \cite{DPModPoly}*{Theorem C}. The intrinsic coefficients are stable under changes of grid, by the comparison
\begin{equation}\label{e:gridchange}
\left|\left\langle f,\Psi^{\tau,1}(R)\right\rangle\right|\lesssim\left|\left\langle f,\Psi^{\tau,1}(P)\right\rangle\right|,\qquad R\subseteq3P,\quad \ell(P)\leq2\ell(R).
\end{equation}
As a consequence, $\mathrm S^\tau_{\mathcal D} f\sim \mathrm S^\tau_{\mathcal D^\star}f$ for every full dyadic grid $\mathcal D^\star$, that is, every image of $\mathcal D$ under a dilation and a translation; for a cube $I$, we write $\mathrm S^\tau_I\coloneqq\mathrm S^\tau_{\mathcal D(I)}$, as $\mathcal D(I)$ is part of such a grid. Write $\tau_1\coloneqq\min\{1,\tau\}$ and $a\coloneqq\lfloor\tau\rfloor$, and let $\omega^\infty_Q$ be the cutoffs of \cite{DPModPoly}*{\S3.1}, equal to one on $\frac32Q$ and supported in $2Q$.

The proof of Theorem \ref{t:square} is modeled after \cite{DPModPoly}*{Section 8}, with Proposition \ref{c:balayage} in place of \cite{DPModPoly}*{Proposition 5.1}. Besides the local weak-type estimate, Proposition \ref{p:wtsquare}, it uses the $L^2$ bound for $\mathrm S^\tau$ and a splitting of the test functions into a local part and tails.
\begin{lemma}\label{l:bessel}
For every $f\in L^2(\R^d)$, $\displaystyle \left\|\mathrm S^\tau  f\right\|_{2}\lesssim\left\|f\right\|_{2}$.
\end{lemma}
\begin{proof}
Choose $\varphi_R\in\Psi^{\tau,1}(R)$ with $|\langle f,\Psi^{\tau,1}(R)\rangle|\leq2|\varphi_R(f)|$, and set $u_R\coloneqq|R|^{1/2}\varphi_R$. Cancellation, regularity and decay give
\[
\left|\int u_R\overline{u_P}\right|
\lesssim
\left(\frac{\min\{\ell(R),\ell(P)\}}{\max\{\ell(R),\ell(P)\}}\right)^{\frac{d+\tau_1}2}
\left(1+\frac{|c(R)-c(P)|}{\max\{\ell(R),\ell(P)\}}\right)^{-d-\frac{\tau_1}2};
\]
see also \cite{FJ90}*{Appendix B}. Apply the Schur test with weights $|R|^{1/2}$: the spatial sum is bounded at each scale, and the remaining scale sum is geometric with exponent $\frac{\tau_1}2$. The resulting Bessel bound gives
\[
\|\mathrm S^\tau f\|_2^2\leq4\sum_R|R||\varphi_R(f)|^2\lesssim\|f\|_2^2.
\]
\end{proof}
\begin{lemma}\label{l:split}
Let $R\subseteq I$ be cubes and $\varphi\in\Psi^{\tau,1}(R)$. There exist $\varphi_j$, $j\geq0$, annihilating $\P^a$, whose partial sums converge to $\varphi$ in $L^2(\R^d)$, such that $\supp\varphi_0\subseteq\frac32I$, $\supp\varphi_j\subseteq2^{j+1}I$, and
\begin{equation}\label{e:split}
\varphi_0\in C\Psi^{\tau,1}(R),\qquad \varphi_j\in C\frac{\ell(R)}{2^{j}\ell(I)}\Phi^{\tau,1}(2^{j+1}I),\quad j\geq1.
\end{equation}
\end{lemma}
\begin{proof}
Let $\Pi^{a+1}_S$ be the adjoint of $\PM^{a+1}_S$, as in \cite{DPModPoly}*{(6.1)}, acting on bounded compactly supported functions $u$. Since $\PM^{a+1}_S$ reproduces $\P^a$, $u-\Pi^{a+1}_Su$ annihilates $\P^a$. Expanding the averaged Taylor formula gives
\begin{equation}\label{e:dualmom}
\Pi^{a+1}_Su=\sum_{|\alpha|\leq a}\ell(S)^{-|\alpha|}
\left(\int(x-c(S))^\alpha u(x)\,\d x\right)\theta_S^\alpha,
\end{equation}
where $\theta_S^\alpha\in C\Phi^{\tau,1}(S)$ are fixed smooth bumps supported in $S$.

Set
\[
\begin{split}
\Sigma_j&\coloneqq\omega^\infty_{\frac34\cdot2^jI}\varphi,
\qquad \Pi_j\coloneqq\Pi^{a+1}_{\frac32\cdot2^jI}\Sigma_j,\\
\varphi_j&\coloneqq(\Sigma_j-\Pi_j)-(\Sigma_{j-1}-\Pi_{j-1}),
\qquad \Sigma_{-1}=\Pi_{-1}\coloneqq0.
\end{split}
\]
The support and moment conditions follow immediately. This is the usual truncation with moment correction, cf.\ \cite{MPTT}*{Lemma 3.1} and \cite{DPOu}*{Lemma 5.2}. We check the estimates with the finite decay available here.

For $j\geq1$, $\Sigma_j-\Sigma_{j-1}$ vanishes on $\frac98\cdot2^{j-1}I\supseteq R$ and is supported in $B=2^{j+1}I$. Rescaling its derivatives and H\"older quotients from $R$ to $B$ costs at most $(\ell(B)/\ell(R))^{d+\tau}$. The decay of $\varphi$ supplies $(\ell(R)/\ell(B))^{d+\tau+1}$, leaving
\begin{equation}\label{e:Sigmaj}
\Sigma_0\in C\Phi^{\tau,1}(R),\qquad
\Sigma_j-\Sigma_{j-1}\in C\frac{\ell(R)}{2^j\ell(I)}\Phi^{\tau,1}(2^{j+1}I),\quad j\geq1.
\end{equation}
Cancellation also gives, for $|\alpha|\leq a$,
\[
\left|\int(x-c(I))^\alpha\Sigma_j(x)\,\d x\right|
=\left|\int(x-c(I))^\alpha(\Sigma_j-\varphi)(x)\,\d x\right|
\lesssim\ell(R)^{\tau+1}(2^j\ell(I))^{|\alpha|-\tau-1}.
\]
Inserting this into \eqref{e:dualmom},
\begin{equation}\label{e:Pij}
\Pi_j\in C\left(\frac{\ell(R)}{2^j\ell(I)}\right)^{\tau+1}
\Phi^{\tau,1}\left(\tfrac32\cdot2^jI\right),\qquad j\geq0.
\end{equation}
For $j=0$ this also places $\Pi_0$ in $C\Phi^{\tau,1}(R)$, since the prefactor supplies the required decay on $\frac32I$. Thus \eqref{e:Sigmaj}--\eqref{e:Pij} imply \eqref{e:split}. Finally, the partial sums $\Sigma_n-\Pi_n$ of the series of the $\varphi_j$ converge to $\varphi$ in $L^2(\R^d)$, by \eqref{e:Pij} and dominated convergence.
\end{proof}
\begin{proposition}\label{p:wtsquare}
For every $0<\eps\leq1$ there exists $C_\eps\lesssim\eps^{-1}$ such that
\[
\left|\left\{x\in Q:\mathrm S^\tau_Q\left(\omega^\infty_Qh\right)(x)>C_\eps\left(\M^\tau_{3Q}h\right)^{*,Q}(r)\right\}\right|\leq\eps|Q|
\]
for all cubes $Q$, all $0<r\leq\frac\eps2$, and all $h\in L^1_{\mathrm{loc}}(\R^d)$ with $h\ind_{3Q}\in L^2(\R^d)$.
\end{proposition}
\begin{proof}
By homogeneity, we can reduce to the case $(\M^\tau_{3Q}h)^{*,Q}(r)=1$.  Indeed, if  $(\M^\tau_{3Q}h)^{*,Q}(r)$ were to vanish instead, so would $\M^\tau_{3Q}h$ at a point of $Q$, whence $h=0$ on $3Q$, by \eqref{e:Mtau} with the cube $3Q$. Thus
\[
U\coloneqq\left\{x\in(3Q)^\circ:\M^\tau_{3Q}h(x)>1\right\}
\]
satisfies $|U\cap Q|\leq r|Q|$.

We use the Calder\'on--Zygmund decomposition of $v\coloneqq\omega^\infty_Qh$ at height one of the proof of \cite{DPModPoly}*{Proposition 8.2}, with Whitney parameter $\Lambda\coloneqq16d$ in place of $8d$:
\[
v=g+\sum_{W\in\mathcal W}b_W,\qquad b_W\coloneqq(v-p_W)\eta^\infty_W,\qquad \|g\|_\infty\lesssim1,\quad\supp g\subseteq3Q,
\]
where the Whitney cubes $W\in\mathcal W$ of $U$ satisfy $\Lambda W\subseteq U$, and $\{\eta^\infty_W\}$ is a smooth partition of unity of $U$ with $\supp\eta^\infty_W\subseteq2W$. Here $p_W$ is the constant for which $b_W$ has vanishing mean, and $b_W=0$ unless $2W\cap2Q\neq\varnothing$, in which case testing \eqref{e:Mtau} at a point of $(3Q)^\circ\setminus U$ close to $W$ gives $|\chi(v)|+|p_W|\lesssim1$ for $\chi\in\Phi^{\tau,1}(W)$ supported in $2W$. Consequently, for a continuous function $\phi$,
\[
\phi(b_W)=\Xi(v)-p_W\int_{\R^d}\Xi,\qquad\Xi\coloneqq\eta^\infty_W\left[\phi-\phi(c(W))\right],
\]
so that, if $\Xi\in C_\phi\Phi^{\tau,1}(W)$,
\begin{equation}\label{e:CZbad}
|\phi(b_W)|\lesssim C_\phi.
\end{equation}
By Lemma \ref{l:bessel} and \eqref{e:gridchange}, $\|\mathrm S^\tau_Qg\|_2\lesssim\|g\|_2\lesssim|Q|^{\frac12}$. For the bad part, let $x\notin U$, so that $D\coloneqq|x-c(W)|\geq8d\ell(W)$, let $R\ni x$ be a cube and $\varphi\in\Psi^{\tau,1}(R)$. The regularity and decay of $\varphi$ on $2W$ give \eqref{e:CZbad} for $\phi=\varphi$, with
\[
C_\varphi\lesssim\left(\frac{\ell(W)}D\right)^{d+\tau_1}\frac{\ell(R)}D\quad\textrm{if }\ell(R)\leq D,\qquad C_\varphi\lesssim\left(\frac{\ell(W)}{\ell(R)}\right)^{d+\tau_1}\quad\textrm{if }\ell(R)>D.
\] Summing the squares over the cubes $R\ni x$ of $\mathcal D(Q)$, one for each sidelength,
\[
\mathrm S^\tau_Q(b_W)(x)\lesssim\left(\frac{\ell(W)}{|x-c(W)|}\right)^{d+\tau_1},\qquad x\in\R^d\setminus U,
\]
and, as $\mathrm S^\tau_Q$ is subadditive,
\[
\|\mathrm S^\tau_Q(v-g)\|_{L^1(Q\setminus U)}\lesssim|U|\lesssim|Q|.
\]
Chebyshev's inequality for $g$ and $v-g$ concludes the proof.
\end{proof}

\subsection{Proof of Theorem \ref{t:square}}\label{ss:pfC}
Fix $0<r\leq r_0$, with $r_0\leq\frac18\eps_d$, and $f\in L^2(\R^d)$ supported in a cube $Q'$. Let $Q_0\coloneqq4Q'$, $\eps\coloneqq\tfrac12\eps_d$ and, for each cube $P$, $h_P\coloneqq f-\PM^m_{3P}f$, $\lambda_P\coloneqq\nat{f}_{m,P}$ and
\[
E_P\coloneqq\left\{x\in P:\M^\tau_{3P}h_P(x)>\lambda_P\right\}\cup\left\{x\in P:\mathrm S^\tau_P\left(\omega^\infty_Ph_P\right)(x)>C_\eps\lambda_P\right\},
\]
with $C_\eps$ as in Proposition \ref{p:wtsquare}, so that $|E_P|\leq(\frac12\eps_d+r)|P|\leq\eps_d|P|$, by Proposition \ref{p:wtsquare} with $h=h_P$, as $r\leq\frac\eps2$, and the definition of $\lambda_P$. For each cube $B$, let
\[
\mu_B\coloneqq\sup\left\{|\varphi(f)|:\;\varphi\in\Psi^{\tau,1}(B),\;\supp\varphi\subseteq B\right\}.
\]
The proof rests on the two estimates
\begin{gather}
\label{e:C1}
\mu_B\lesssim\lambda_P,\qquad B\subseteq3P,\quad 3B\cap(P\setminus E_P)\neq\varnothing,
\\
\label{e:C2}
\sum_{R\in\mathcal D(Q),\,R\ni y}\left|\left\langle f,\Psi^{\tau,1}(R)\right\rangle\right|^2\lesssim\lambda_Q^2+\sum_{j\geq1}2^{-j}\mu_{2^jQ}^2,\qquad y\in Q\setminus E_Q,
\end{gather}
for $P\in\mathcal D^+(Q_0)$ and $Q\in\mathcal D(Q_0)$, respectively.

We verify \eqref{e:C1}. Let $y\in3B\cap(P\setminus E_P)$. By \cite{DPModPoly}*{(4.3) and (4.1)}, there is a cube $B'$ with $B\cup\{y\}\subseteq B'\subseteq3P$ and $\ell(B')\leq3\ell(B)$, so that each $\varphi$ in the definition of $\mu_B$ is a bounded multiple of a test function of \eqref{e:Mtau} relative to $B'$. As $\varphi$ annihilates $\P^a\supseteq\P^{m-1}$ and $f-h_P\in\P^{m-1}$,
\[
|\varphi(f)|=|\varphi(h_P)|\lesssim\M^\tau_{3P}h_P(y)\leq\lambda_P.
\]
We verify \eqref{e:C2}. Let $R\in\mathcal D(Q)$ contain $y$, and split $\varphi\in\Psi^{\tau,1}(R)$ as in Lemma \ref{l:split}, with $I=Q$. As $\varphi_0$ annihilates $\P^{m-1}$ and is supported in $\frac32Q$, where $\omega^\infty_Q=1$, while each $\varphi_j$, $j\geq1$, is a multiple of a test function in the definition of $\mu_{2^{j+1}Q}$,
\[
|\varphi_0(f)|=\left|\varphi_0\left(\omega^\infty_Qh_Q\right)\right|\lesssim\left|\left\langle\omega^\infty_Qh_Q,\Psi^{\tau,1}(R)\right\rangle\right|,\qquad |\varphi_j(f)|\lesssim\frac{\ell(R)}{2^j\ell(Q)}\mu_{2^{j+1}Q}.
\]
Summing the squares over the cubes $R\ni y$ of $\mathcal D(Q)$, one for each sidelength, and using the Cauchy--Schwarz inequality,
\[
\sum_{R\in\mathcal D(Q),\,R\ni y}\left|\left\langle f,\Psi^{\tau,1}(R)\right\rangle\right|^2\lesssim\left[\mathrm S^\tau_Q\left(\omega^\infty_Qh_Q\right)(y)\right]^2+\sum_{j\geq1}2^{-j}\mu^2_{2^{j+1}Q},
\]
which is \eqref{e:C2}, as $y\notin E_Q$, after a shift of the index in the last sum.

\emph{Stopping.} By \eqref{e:C1} and \eqref{e:C2}, the assumption \eqref{e:P4} of Proposition \ref{c:balayage} holds for $a_R\coloneqq|\langle f,\Psi^{\tau,1}(R)\rangle|^2$, with $\lambda_P^2$ in place of $\lambda_P$. Indeed, if $P\in\mathcal T_k(Q)$, then $2^kQ\subseteq3P$, and \eqref{e:C1} with $B=2^kQ$ yields $\mu_{2^kQ}\lesssim\lambda_P$. The tailed local bounds are finite, as $\lambda_P\lesssim r^{-1}\langle|f|\rangle_{3P}$ by \eqref{e:poelem}, which is bounded uniformly over the cubes $P=2^jQ_0$. Therefore, Proposition \ref{c:balayage} yields a $\frac12$-sparse collection $\mathcal S\subseteq\mathcal D^+(Q_0)$, containing the cubes $2^jQ_0$, $j\geq1$, such that
\begin{equation}\label{e:C3}
\sum_{R\in\mathcal D(Q_0)}\left|\left\langle f,\Psi^{\tau,1}(R)\right\rangle\right|^2\ind_R\lesssim\sum_{S\in\mathcal S}\lambda_S^2\ind_S\qquad\textrm{a.e.\ on } \R^d.
\end{equation}

\emph{Exterior cubes.} It remains to estimate the cubes outside $\mathcal D(Q_0)$. By the comparison following \eqref{e:gridchange}, $\mathcal D$ may be replaced by a dyadic grid $\mathcal D^\star\supseteq\mathcal D(Q_0)$. Fix $x\in\R^d$, and let $R\in\mathcal D^\star\setminus\mathcal D(Q_0)$ contain $x$, so that $R$ either strictly contains $Q_0$ or does not meet it. With $\rho\coloneqq\dist_\infty(x,Q')$, where $\dist_\infty$ is distance in the maximum norm, split $\varphi\in\Psi^{\tau,1}(R)$ as in Lemma \ref{l:split}, where
\[
I\coloneqq\begin{cases}\textrm{the smallest of the cubes } 2^lQ_0 \textrm{ containing } R, & \ell(R)\geq\frac14\rho,\\
\textrm{the cube of center } x \textrm{ and sidelength } \frac12\rho, & \ell(R)<\frac14\rho.\end{cases}
\]
In the first case $\ell(I)\sim\ell(R)$; in the second, $x\notin Q_0$ and $\frac32I$ does not meet $Q'$, whence $\varphi_0(f)=0$. For $j\geq0$, let $P_j$ be the smallest of the cubes $2^lQ_0$, $l\geq1$, containing $2^{j+1}I$; it contains $x$, and $\ell(P_j)\sim2^j\ell(I)$. Then, for $j\geq1$, and for $j=0$ in the first case,
\[
\varphi_j\in C\frac{\ell(R)}{2^j\ell(I)}\Psi^{\tau,1}(P_j),\qquad\supp\varphi_j\subseteq P_j,
\]
and \eqref{e:C1} with $B=P=P_j$ yields
\begin{equation}\label{e:C4}
\left|\left\langle f,\Psi^{\tau,1}(R)\right\rangle\right|\lesssim\frac{\ell(R)}{\ell(I)}\sum_{j\geq0}2^{-j}\lambda_{P_j}.
\end{equation}
In the first case $P_j=2^{j+1}I$, and boundedly many cubes $R$ share the same $I$; in the second, $I$ does not depend on $R$, and the ratios $\ell(R)/\ell(I)<\frac12$ form a geometric sequence, there being one cube $R\ni x$ of each sidelength $\ell(R)<\frac14\rho$. Squaring \eqref{e:C4}, applying the Cauchy--Schwarz inequality and summing over $R$,
\[
\sum_{R\in\mathcal D^\star\setminus\mathcal D(Q_0),\,R\ni x}\left|\left\langle f,\Psi^{\tau,1}(R)\right\rangle\right|^2\lesssim\sum_{l\geq1}\lambda_{2^lQ_0}^2\ind_{2^lQ_0}(x),
\]
which, together with \eqref{e:C3}, gives the first estimate of Theorem \ref{t:square}. The second holds for the same collection. Indeed, by the proof of Proposition \ref{c:balayage}, the cubes of $\mathcal S$ in $\mathcal D(Q_0)$ form the collection of \cite{DPModPoly}*{Proposition 5.1} for the sets $E_Q$, $Q\in\mathcal D(Q_0)$, which depends on these sets alone. For each such $Q$, the set $E_Q$ contains the exceptional set $\{x\in Q:\M^\tau_{3Q}h_Q(x)>\lambda_Q\}$ of the proof of \cite{DPModPoly}*{Proposition 9.5}, which uses the same $h_Q$, $\lambda_Q$ and $Q_0=4Q'$, and $|E_Q|\leq\eps_d|Q|$. By \cite{DPModPoly}*{Remark 5.2}, that proof applies with the sets $E_Q$ in place of its own, and yields the collection $\mathcal S$.

\begin{remark}\label{r:square}
Theorem \ref{t:square} only uses the vanishing of the moments of order less than $\max\{1,m\}$: the vanishing mean in Lemma \ref{l:bessel}, and the vanishing of the moments of order less than $m$ in \eqref{e:C1}. Lemma \ref{l:split} holds with $\P^a$ replaced by $\P^{\max\{0,m-1\}}$, and $\Psi^{\tau,1}(R)$ by the class of those $\varphi\in\Phi^{\tau,1}(R)$ annihilating $\P^{\max\{0,m-1\}}$, as long as the moments involved converge, that is, if $\tau>m-2$. Therefore, $\Psi^{\tau,1}$ may be replaced by the latter class, both in \eqref{e:intrinsic} and in the definition of $\mu_B$; it contains $\Psi^{\tau,1}$ when $\tau>m-1$, and the restriction $\tau>m-1$ may then be relaxed to $\tau>m-2$.
\end{remark}

\section{Potential and weighted bounds}\label{s:cons}
This section contains two weighted lemmas, the proofs of Corollaries \ref{c:sqpot} and \ref{c:weighted}, and the sharpness of the latter. The first lemma is the quadratic counterpart of \cite{DPModPoly}*{Lemma 9.3}.
\begin{lemma}\label{l:sqAinf}
Let $0<p<\infty$, $w\in A_\infty$, and let $\mathcal S$ be an $\eta$-sparse collection of cubes. If $\lambda_S\geq0$ for $S\in\mathcal S$, then
\[
\Big\|\Big(\sum_{S\in\mathcal S}\lambda_S^2\ind_S\Big)^{\frac12}\Big\|_{L^p(w)}\lesssim_p[w]_{A_\infty}^{\frac12}\Big\|\sup_{S\in\mathcal S}\lambda_S\ind_S\Big\|_{L^p(w)}.
\]
The implied constants also depend on $\eta$, and grow at most like $\sqrt p$ as $p\to\infty$.
\end{lemma}
\begin{proof}
Apply \cite{DPModPoly}*{Lemma 9.3} to the coefficients $\lambda_S^2$ at exponent $p/2$, and take square roots:
\[
\left\|\sum_{S\in\mathcal S}\lambda_S^2\ind_S\right\|_{L^{p/2}(w)}^{1/2}
\lesssim_p[w]_{A_\infty}^{1/2}
\left\|\sup_{S\in\mathcal S}\lambda_S^2\ind_S\right\|_{L^{p/2}(w)}^{1/2}.
\]
Its constant is $O(p)$ before taking square roots, giving the asserted $O(\sqrt p)$ growth. Infinite collections follow by monotone convergence.
\end{proof}
The second lemma makes the comparison of \cite{DPModPoly}*{\S2.7} between the maximal percentile oscillations and the smooth grand maximal function explicit in the weight.
\begin{lemma}\label{l:median}
Let $m\geq0$ be an integer, $0<p<\infty$ and $w\in A_q$, $1<q<\infty$. For every $f\in L^1_{\mathrm{loc}}(\R^d)$,
\[
\left\|\mathcal N_mf\right\|_{L^p(w)}\lesssim\left(\frac Cr\right)^{\frac qp}[w]_{A_q}^{\frac1p}\big\|\M^\tau_{\R^d}f\big\|_{L^p(w)},
\]
where $C$ depends only on $d$.
\end{lemma}
\begin{proof}
By \eqref{e:poelem} with $\pi=0$, and as $\M^\tau_{3Q}f\leq\M^\tau_{\R^d}f$ on $Q$, by \eqref{e:Mtau}, for every cube $Q$,
\[
\nat{f}_{m,Q}\lesssim\nat{f}_{0,Q}\leq g^{*,Q}(r),\qquad g\coloneqq\M^\tau_{\R^d}f.
\]
For a measurable set $E$, let $U_E$ be the union of the cubes $Q$ with $|E\cap Q|>r|Q|$, which is open, as each such $Q$ may be slightly dilated. As $g^{*,Q}(r)>t$ if and only if $|\{g>t\}\cap Q|>r|Q|$, it suffices to show that
\[
w(U_E)\leq\left(\frac{9^d}r\right)^q[w]_{A_q}\,w(E),
\]
and to integrate in $t>0$ against $pt^{p-1}$, with $E=\{g>t\}$. By monotone convergence, we may assume that $E$ is bounded. By the three lattice theorem \cite{LN19}*{Theorem 3.1 and Remark 3.2}, there are $3^d$ dyadic lattices such that each cube $Q$ is contained in a cube $P$ of one of them with $|P|\leq3^d|Q|$, whence $|E\cap P|>3^{-d}r|P|$ if $|E\cap Q|>r|Q|$. As $E$ is bounded, the cubes $P$ of each lattice with the latter property are contained in maximal ones, which are pairwise disjoint and, by H\"older's inequality and the definition of $[w]_{A_q}$, satisfy
\[
w(P)\leq[w]_{A_q}\left(\frac{|P|}{|E\cap P|}\right)^qw(E\cap P)\leq[w]_{A_q}\left(\frac{3^d}r\right)^qw(E\cap P).
\]
Summing over the maximal cubes, and then over the $3^d$ lattices, gives the claim, as $3^d\leq3^{dq}$.
\end{proof}
In the rest of this section, $\mathcal S$ is the $\frac12$-sparse collection of Theorem \ref{t:square} constructed in \S\ref{ss:pfC}; as $\mathcal S\subseteq\mathcal D^+(Q_0)$, its cubes are pairwise nested or disjoint, and each point of $\R^d$ belongs to at most one of them with sidelength in any given interval $[2^k,2^{k+1})$.

\begin{proof}[Proof of Corollary \ref{c:sqpot}]
The bound \eqref{e:sqfrac} holds with the collection $\mathcal S$ of Theorem \ref{t:square}, by Theorem \ref{t:square} and \cite{DPModPoly}*{(6.7)}, which also gives, together with Jensen's inequality,
\[
\begin{split}
\nat{f}_{m,S}&\lesssim\frac{\ell(S)^m}r\left\langle|\nabla^mf|\right\rangle_{3S},\\
\left[\nat{f}_{m,S}\right]^2&\lesssim\frac{\ell(S)^{2m}}{r^2}\left\langle|\nabla^mf|^2\right\rangle_{3S},\qquad S\in\mathcal S.
\end{split}
\]
The first potential bound then follows from Theorem \ref{t:square}, the embedding $\ell^1\subseteq\ell^2$, and \cite{DPModPoly}*{Lemma 4.1} with $s=m$, and the second from \cite{DPModPoly}*{Lemma 4.1} with $s=2m$ and $|\nabla^mf|^2$ in place of $g$.
\end{proof}

\begin{proof}[Proof of Corollary \ref{c:weighted}]
Let first $f$ have compact support. As $\nat{f}_{m,S}\ind_S\leq\mathcal N_mf$, Theorem \ref{t:square} and Lemma \ref{l:sqAinf}, with $\lambda_S\coloneqq\nat{f}_{m,S}$, yield \eqref{e:weighted}, with the growth $\sqrt p$ of the constant.

For general $f\in L^2$, fix a cube $Q_0$. The local construction in \S\ref{ss:pfC} still gives \eqref{e:C3}, since
$\lambda_P\lesssim r^{-1}|3P|^{-1/2}\|f\|_2$ makes the tails finite. On $Q_0$, the remaining cubes of a grid containing $\mathcal D(Q_0)$ are ancestors of $Q_0$. Apply Lemma \ref{l:split} with $I$ the least concentric dilate $2^lQ_0$ containing each such ancestor; then $\ell(I)\sim\ell(R)$ and \eqref{e:C1} gives \eqref{e:C4}. The same geometric summation proves the square bound on $Q_0$. Lemma \ref{l:sqAinf} therefore bounds $\|\ind_{Q_0}\mathrm S^\tau f\|_{L^p(w)}$ uniformly in $Q_0$. Let $Q_0$ increase to $\R^d$. The last claim follows, as $\mathcal N_mf\lesssim r^{-1}\M_m(|\nabla^mf|)$ \cite{DPModPoly}*{Remark 2.5}.
\end{proof}
\begin{remark}\label{r:sharp}
The factor $[w]_{A_\infty}^{1/2}$ and the growth $\sqrt p$ in Corollary \ref{c:weighted} cannot be improved. Let $f\coloneqq\ind_{[0,1)^d}$ and $R_j\coloneqq[0,2^{-j})^d$. Fix $\xi\in\mathcal{C}^\infty_c((-\frac34,\frac34))$ with $\xi^{(a)}(-\frac12)\neq0$, so small that the function
\[
\varphi(x)\coloneqq\xi^{(a+1)}(x_1)\cdots\xi^{(a+1)}(x_d),
\]
whose moments of order at most $a$ vanish, satisfies $\varphi_R\in\Psi^{\tau,1}(R)$ for every cube $R$. For $j\geq1$, as $2^j-\frac12>\frac34$,
\[
\varphi_{R_j}(f)=\Big(\int_{-1/2}^{2^j-1/2}\xi^{(a+1)}\Big)^d=\left(-\xi^{(a)}(-\tfrac12)\right)^d\neq0,
\]
so that
\[
\mathrm S^\tau f(x)^2\gtrsim\#\{j\geq1:x\in R_j\}\gtrsim\log\frac1{|x|},\qquad x\in(0,\tfrac14)^d,
\]
while $\mathcal N_mf\lesssim r^{-1}\M f\lesssim r^{-1}\min\{1,|x|^{-d}\}$, by \eqref{e:poelem}. The power weights $w_\varepsilon(x)\coloneqq|x|^{-d+\varepsilon}$, $0<\varepsilon\leq\frac12$, satisfy $[w_\varepsilon]_{A_\infty}\sim\varepsilon^{-1}$ \cite{DPModPoly}*{Lemma 11.5}. For $0<p<\infty$ and $\varepsilon\leq\frac{dp}2$, polar coordinates and the change of variables $t=e^{-s/\varepsilon}$ give, with implicit constants depending on $p$ only if $p<2$,
\begin{align*}
\left\|\mathcal N_mf\right\|_{L^p(w_\varepsilon)}&\lesssim\frac1r\,\varepsilon^{-\frac1p},\\
\left\|\mathrm S^\tau f\right\|_{L^p(w_\varepsilon)}&\gtrsim\Big(\int_0^{1/4}\Big(\log\frac1t\Big)^{\frac p2}t^{\varepsilon-1}\,\d t\Big)^{\frac1p}\gtrsim\varepsilon^{-\frac1p-\frac12}\,\Gamma\big(1+\tfrac p2\big)^{\frac1p},
\end{align*}
so that the ratio of the two norms is $\gtrsim_p r\,\varepsilon^{-1/2}\sim r\,[w_\varepsilon]_{A_\infty}^{1/2}$, and $\gtrsim r\sqrt p\,[w_\varepsilon]_{A_\infty}^{1/2}$ for $p\geq2$.
\end{remark}

\section{Proof of Theorem \ref{t:SoT}}\label{s:SoT}
The proof of Theorem \ref{t:SoT} rests on Corollary \ref{c:weighted} and on a pointwise inequality between maximal percentile oscillations, Proposition \ref{p:SoT}, for the model wavelet operators of \cite{DPModPoly}*{Definition 7.2}: given packets $\varphi_R\in\Phi^{\sigma,\rho}(R)$ and $\psi_R\in\Phi^{\sigma',\rho'}(R)$, $R\in\mathcal D$, with $\sigma,\sigma'>0$ and $\rho,\rho'\geq0$, the operator
\[
Vf\coloneqq\sum_{R\in\mathcal D}|R|\varphi_R(f)\psi_R
\]
is a \emph{$(\sigma,\sigma';\rho,\rho')$-model} if its partial sums over the finite subcollections of $\mathcal D$ have norm at most one on $L^2(\R^d)$; the $\varphi_R$ are its \emph{input packets}, the $\psi_R$ its \emph{output packets}, and $K$ is its kernel \cite{DPModPoly}*{\S7.7}. In this section,
\[
\tau<\min\{1,\sigma,\sigma'+\rho'\},\qquad\nu\coloneqq\tfrac12\min\{1,\sigma',\sigma+\rho-\tau,\sigma'+\rho'-\tau\};
\]
the condition on $\tau$ implies the range condition \cite{DPModPoly}*{(7.10)}. In the next lemma, the vanishing mean of the output packets plays the role of the condition $T^*1=0$.
\begin{lemma}\label{l:adjoint}
Let $V$ be a model wavelet operator whose output packets have vanishing mean. For every cube ${R_0}$ and every $\zeta\in\Phi^{\tau,1}({R_0})$ with $\supp\zeta\subseteq{R_0}$, the function
\[
V^*\zeta\coloneqq\sum_{R\in\mathcal D}|R|\psi_R(\zeta)\varphi_R
\]
satisfies, for all $z,z'\in\R^d$,
\begin{equation}\label{e:adjoint}
\left|V^*\zeta(z)\right|\lesssim\frac1{|{R_0}|},\qquad \left|V^*\zeta(z)-V^*\zeta(z')\right|\lesssim\frac1{|{R_0}|}\left(\frac{|z-z'|}{\ell({R_0})}\right)^\tau.
\end{equation}
\end{lemma}
\begin{proof}
Let $s\coloneqq\ell({R_0})$, $\sigma_1\coloneqq\min\{1,\sigma\}$, and $\ell(R)=q$. The size of $\psi_R$ if $q\geq s$, and, if $q<s$, its vanishing mean and decay, as $\tau<\sigma'+\rho'$, and the H\"older continuity of $\zeta$, give
\[
|\psi_R(\zeta)|\lesssim s^{-d}\min\left\{\left(\frac qs\right)^\tau,\left(\frac sq\right)^d\right\}.
\]
As the centers of the cubes of sidelength $q$ form a translate of $q\Z^d$, the input packets satisfy, for $|z-z'|\leq q$,
\[
\sum_{\ell(R)=q}|R|\left|\varphi_R(z)\right|\lesssim1,\qquad \sum_{\ell(R)=q}|R|\left|\varphi_R(z)-\varphi_R(z')\right|\lesssim\left(\frac{|z-z'|}q\right)^{\sigma_1}.
\]
Summing the absolute values over $q$ gives \eqref{e:adjoint}. For the second estimate, which follows from the first if $|z-z'|>s$, let $|z-z'|\leq s$; then the scales $q<|z-z'|$ are estimated by the first lattice sum at $z$ and at $z'$, and the others by the second, as $\sigma_1>\tau$.
\end{proof}
\begin{proposition}\label{p:SoT}
Let $V$ be a model wavelet operator whose output packets have vanishing mean, and let $m_0\in\{0,1\}$. If $m_0=1$, assume that the input packets have vanishing mean as well. For every $f\in L^2(\R^d)$ and every cube $Q$, there is a constant $c_Q$ such that, with $r'\coloneqq2^{-2d-1}r$,
\begin{equation}\label{e:SoTperc}
\left(\M^\tau_{3Q}\left[Vf-c_Q\right]\right)^{*,Q}(r)\lesssim\frac1r\sum_{j\geq0}2^{-j\nu}\natl{f}{r'}_{m_0,2^{j+2}Q}.
\end{equation}
In particular, $\mathcal N^{[r]}_1(Vf)\lesssim r^{-1}\mathcal N^{[r']}_{m_0}f$.
\end{proposition}
\begin{proof}
Set $\widetilde Q=4Q$, $\varpi=\omega^\infty_{\widetilde Q}$, $c=\PM^{m_0}_{3\widetilde Q}f$ and $h=f-c$. Thus $c$ is constant, and zero if $m_0=0$. Decompose
\[
Vf=V(\varpi h)+G,\qquad G=V(\varpi c)+V((1-\varpi)f).
\]
We estimate the first term weakly and the oscillation of the second uniformly.

\emph{Local term.} With $\eps=4^{-d}r=2r'$, we claim that
\begin{equation}\label{e:SoTwt}
\left|\left\{x\in\widetilde Q:
\M^\tau_{3\widetilde Q}[V(\varpi h)](x)
>C_\eps\natl{f}{r'}_{m_0,\widetilde Q}\right\}\right|
\leq\eps|\widetilde Q|,\qquad C_\eps\lesssim\eps^{-1}.
\end{equation}
Use the decomposition in Proposition \ref{p:wtsquare} for $v=\varpi h$, after normalizing $\natl{f}{r'}_{m_0,\widetilde Q}=1$. Its good part is bounded in $L^2$ by $C|\widetilde Q|^{1/2}$, and both $V$ and $\M^\tau_{3\widetilde Q}\leq\M$ are $L^2$ bounded. For the bad part it suffices to prove
\begin{equation}\label{e:SoTbad}
\M^\tau_{3\widetilde Q}(Vb_W)(x)
\lesssim\left(\frac{\ell(W)}{|x-c(W)|}\right)^{d+\tau},
\qquad x\in3\widetilde Q\setminus U.
\end{equation}
Write $D=|x-c(W)|\geq8d\ell(W)$, and let $\zeta$ be a test supported in $R_0\ni x$.
If $\sqrt d\ell(R_0)\leq D/4$, then $R_0$ and $2W$ are separated. The kernel representation and the kernel estimates \cite{DPModPoly}*{(7.16) and (7.11)--(7.12)} place
$\eta_W^\infty[K(x',\cdot)-K(x',c(W))]$ in
$C(\ell(W)/D)^{d+\tau}\Phi^{\tau,1}(W)$, uniformly for $x'\in R_0$. Apply \eqref{e:CZbad} and integrate against $\zeta$.
Otherwise, use $\zeta(Vb_W)=b_W(V^*\zeta)$. Lemma \ref{l:adjoint} places
\[
\eta_W^\infty[V^*\zeta-V^*\zeta(c(W))]
\in C\frac{|W|}{|R_0|}
\left(\frac{\ell(W)}{\ell(R_0)}\right)^\tau\Phi^{\tau,1}(W),
\]
and \eqref{e:CZbad} again gives \eqref{e:SoTbad}. Its integral off $U$, summed over $W$, is $\lesssim|\widetilde Q|$. Chebyshev's inequality proves \eqref{e:SoTwt}.

\emph{Far term.} We prove
\begin{equation}\label{e:SoTfar}
|G(x')-G(y)|\lesssim
\sum_{j\geq0}2^{-j\nu}\natl{f}{r'}_{m_0,2^j\widetilde Q},
\qquad\textrm{a.e. }x',y\in3Q.
\end{equation}
First,
\begin{equation}\label{e:SoTG}
G(x')-G(y)=\int[K(x',z)-K(y,z)](1-\varpi(z))(f(z)-c)\,\d z.
\end{equation}
The integral converges absolutely by \cite{DPModPoly}*{(7.11), (7.13)}. When $m_0=0$, this is the kernel representation. When $m_0=1$, it follows by cutting off the constant and using the cancellation of the input packets. Indeed, with $\iota=\frac12\min\{1,\sigma+\rho\}$, their pairing with $\omega^\infty_{AQ}$ is bounded by
$C\min\{(q/(A\ell(Q)))^\iota,(A\ell(Q)/q)^d\}$ at scale $q$.
The lattice sum of the output differences on $3Q$ is bounded by
$C\min\{1,(\ell(Q)/q)^\nu\}$. Summing their product gives
\[
\sup_{x',y\in3Q}
|V(\omega^\infty_{AQ})(x')-V(\omega^\infty_{AQ})(y)|
\lesssim(1+\log A)A^{-\min\{\iota,\nu\}}\longrightarrow0.
\]
Apply the kernel representation to $\omega^\infty_{AQ}-\varpi$ and let $A\to\infty$ to obtain \eqref{e:SoTG}.

Write
\[
1-\varpi=\sum_{j\geq0}u_j,\qquad u_j\coloneqq\omega^\infty_{2^{j+1}\widetilde Q}-\omega^\infty_{2^j\widetilde Q},\qquad P_j\coloneqq2^{j+2}\widetilde Q.
\]
The support of $u_j$ lies in $P_j\setminus\frac32\cdot2^j\widetilde Q$. The size, regularity and mixed kernel estimates \cite{DPModPoly}*{(7.11)--(7.14)} give
\[
u_j[K(x',\cdot)-K(y,\cdot)]
\in C2^{-j\nu}\Phi^{\tau,1}(P_j),\qquad 2^j\gg1.
\]
For the boundedly many small $j$, use the size and input regularity estimates separately at $x'$ and $y$. Testing at every point of $P_j$ in \eqref{e:SoTG},
\[
|G(x')-G(y)|\lesssim
\sum_{j\geq0}2^{-j\nu}\inf_{P_j}\M^\tau_{3P_j}(f-c).
\]
For $m_0=0$, each infimum is at most $\natl{f}{r'}_{0,P_j}$.
For $m_0=1$, put $c_i=\PM^1_{3\cdot2^i\widetilde Q}f$. By \cite{DPModPoly}*{(4.6) and Lemma 4.3(ii)},
\[
|c_i-c_{i+1}|\lesssim\natl{f}{r'}_{1,2^{i+1}\widetilde Q},\qquad
\inf_{P_j}\M^\tau_{3P_j}(f-c)
\lesssim\sum_{i=0}^{j+2}\natl{f}{r'}_{1,2^i\widetilde Q}.
\]
Summing the geometric series proves \eqref{e:SoTfar} in both cases.

Choose $y\in3Q$ for which \eqref{e:SoTfar} holds almost everywhere in $x'$, and set $c_Q=G(y)$. As $\M^\tau_{3Q}\leq\M^\tau_{3\widetilde Q}$ on $Q$, and the test functions in \eqref{e:Mtau} have $L^1$ norm at most one, \eqref{e:SoTwt}--\eqref{e:SoTfar}, with $\eps|\widetilde Q|=r|Q|$, prove \eqref{e:SoTperc}. Finally use \eqref{e:poelem} and take the supremum over $Q$.
\end{proof}
\begin{proof}[Proof of Theorem \ref{t:SoT}]
Let $s$ be an admissible level at which Corollary \ref{c:weighted} holds with $m=1$, let $n$ be twice the number of multi-indices of order $\kappa$, and let $f\in \mathcal{C}^\infty_c(\R^d)$ and $0<r\leq r_0\leq s'$, where $s'$ is the level $r'$ of Proposition \ref{p:SoT} relative to $s/n$.

\emph{Sufficiency.} Fix $\tau<\delta_1<\delta$, and let $V_{\alpha,\beta}$ be the cancellative model sums of \cite{DPModPoly}*{Proposition 7.6}, with input packets $\varphi^{\alpha,\beta}_R\in C\Psi^{\delta_1,0}(R)$ and output packets $\psi_R$. Under the conditions in (a), by \cite{DPModPoly}*{Proposition 7.6 and (7.8)}, the lower paraproduct forms vanish, the top ones are given by the canonical paraproducts with symbols $\mathfrak b^{\alpha,\gamma}\in\mathrm{BMO}(\R^d)$, and, if $\kappa=0$, the adjoint term vanishes. Thus
\[
\partial^\alpha Tf=\sum_{|\beta|=\kappa}V_{\alpha,\beta}\left(\partial^\beta f\right)+\sum_{|\gamma|=\kappa}\frac1{\gamma!}P_{\mathfrak b^{\alpha,\gamma}}\left(\partial^\gamma f\right),
\]
all terms belonging to $L^2(\R^d)$. Up to constant multiples, admissible for $V_{\alpha,\beta}$ and $\lesssim\|b^{\alpha,\gamma}\|_{\mathrm{BMO}}$ for $P_{\mathfrak b^{\alpha,\gamma}}$, each of the $n$ operators on the right hand side is a model whose output packets, the wavelets $\psi_R$, have vanishing mean, and for which $\tau$ is in the range of Proposition \ref{p:SoT}, by \cite{DPModPoly}*{Remark 7.4 and Lemma 10.7(i)}. As the percentile of a sum of $n$ nonnegative functions at the level $s$ is at most the sum of their percentiles at the level $s/n$, Proposition \ref{p:SoT} with $m_0=0$, at the level $s/n$, and the monotonicity of the percentiles in the level give
\begin{equation}\label{e:SoTpoint}
\mathcal N^{[s]}_1\left(\partial^\alpha Tf\right)\lesssim\sum_{|\beta|=\kappa}\mathcal N^{[s']}_0\left(\partial^\beta f\right)\leq\sum_{|\beta|=\kappa}\mathcal N^{[r]}_0\left(\partial^\beta f\right).
\end{equation}
Corollary \ref{c:weighted}, with $m=1$ and at the level $s$, applied to $\partial^\alpha Tf$, and \eqref{e:SoTpoint}, whose implied constant is thus $\lesssim1+\max_{|\gamma|=\kappa}\|b^{\alpha,\gamma}\|_{\mathrm{BMO}}$, give \eqref{e:SoT}, with constants growing at most like $\sqrt p$, the bound being trivial when its right hand side is infinite. Under the conditions in (b), the paraproducts are absent, the input packets of $V_{\alpha,\beta}$ have vanishing mean, and Proposition \ref{p:SoT} with $m_0=1$ gives \eqref{e:SoTpoint} with $\mathcal N_1(\partial^\beta f)$ on its right hand side, whence \eqref{e:SoT1}. The factor $[w]_{A_\infty}^{1/2}$ cannot be improved for the identity, which satisfies the conditions in (a) and (b) with $\kappa=0$, by Remark \ref{r:sharp}, applied to smooth approximations of $\ind_{[0,1)^d}$.

\emph{Necessary lower and top conditions.} Assume \eqref{e:SoT}. The wavelets $\psi_S$ belong to $C\Psi^{\tau,1}(S)$, so that, for $g\in L^2(\R^d)$ and $h\in \mathcal{C}^\infty_c(\R^d)$, by \eqref{e:onb},
\[
|\langle g,h\rangle|\leq\sum_{S\in\mathcal D}|S||\psi_S(g)||\psi_S(h)|\lesssim\int_{\R^d}\mathrm S^\tau g\,\mathrm S^\tau h.
\]
As $\|\mathrm S^\tau h\|_{p'}\lesssim_p\|h\|_{p'}$ for $1<p<\infty$, by Corollary \ref{c:weighted} with $m=0$, $w\equiv1$ and at a fixed admissible level, and \eqref{e:poelem}, it follows that $\|g\|_p\lesssim_p\|\mathrm S^\tau g\|_p$. With $g=\partial^\alpha Tf$, \eqref{e:SoT} with $w\equiv1$ and $\mathcal N_0(\partial^\beta f)\lesssim r^{-1}\M(\partial^\beta f)$ give $\|\partial^\alpha Tf\|_p\lesssim_p r^{-1}\|\nabla^\kappa f\|_p$ for $1<p<\infty$. By \cite{DPModPoly}*{Lemma 11.2}, the lower symbols vanish and the top symbols belong to $\mathrm{BMO}$; if $\kappa=0$, so does $b^{0,\star}$. 

\emph{The adjoint symbol.} Suppose $\kappa=0$. By \cite{DPModPoly}*{Proposition 7.6},
\[
Tf=V_{0,0}f+P_{\mathfrak b}f+P^*_{\mathfrak b^\star}f,\qquad\mathfrak b\coloneqq\mathfrak b^{0,0},\quad\mathfrak b^\star\coloneqq\mathfrak b^{0,\star}\in\mathrm{BMO}(\R^d).
\] Let in addition $\int f=0$. Then $\mathcal N_0f$ is integrable, by \cite{DPModPoly}*{Proposition 2.8}, as $f\in H^1(\R^d)$, hence so are $\mathrm S^\tau(Tf)$, by \eqref{e:SoT} with $p=1$ and $w\equiv1$, $\mathrm S^\tau(V_{0,0}f)$ and $\mathrm S^\tau(P_{\mathfrak b}f)$, by the proof of the sufficiency, and $\mathrm S^\tau(P^*_{\mathfrak b^\star}f)$, by subadditivity. As the numbers $|Q||\psi_Q(f)|$ are summable over $Q\in\mathcal D$, by the vanishing mean of $\psi_Q$ at small scales and of $f$ at large ones, the series \cite{DPModPoly}*{(7.8)} defining $P^*_{\mathfrak b^\star}f$ converges absolutely in $L^1(\R^d)$. Testing \eqref{e:intrinsic} with rescalings of a fixed function with vanishing mean, equal to a nonzero constant on a large ball, at the cubes containing $x$ of sidelength at least $16|x|$, shows that $\mathrm S^\tau g(x)\gtrsim|\int g||x|^{-d}$ for $|x|$ large, depending on $g$, if $g\in L^1(\R^d)\cap L^2(\R^d)$. Therefore, as $\int\vartheta_Q=1$,
\[
0=\int_{\R^d}P^*_{\mathfrak b^\star}f=\sum_{Q\in\mathcal D}|Q|\psi_Q(\mathfrak b^\star)\psi_Q(f).
\]
Applied to mollifications of $\psi_S$, which have vanishing mean, are supported in a fixed cube and are uniformly Lipschitz, so that their coefficients are dominated uniformly, this gives $b^{0,\star}_S=\psi_S(\mathfrak b^\star)=0$ in the limit, as $\psi_Q(\psi_S)=\delta_{QS}|S|^{-1}$ by orthonormality.

\emph{Vanishing of the top symbols in (b).} Assume \eqref{e:SoT1}, with constants $C_p$. As $\mathcal N_1g\lesssim\mathcal N_0g$ by \eqref{e:poelem}, \eqref{e:SoT} holds as well, so that the conditions in (a) hold, by their necessity, and it remains to show that the top symbols vanish. By the representation of $\partial^\alpha Tf$ in the proof of the sufficiency and \eqref{e:onb}, for every $R\in\mathcal D$,
\begin{equation}\label{e:SoTcoef}
\psi_R\left(\partial^\alpha Tf\right)=\sum_{|\beta|=\kappa}\varphi^{\alpha,\beta}_R\left(\partial^\beta f\right)+\sum_{|\gamma'|=\kappa}\frac1{\gamma'!}\,b^{\alpha,\gamma'}_R\,\vartheta_R\left(\partial^{\gamma'}f\right).
\end{equation}
Fix $R\in\mathcal D$, $|\gamma|=\kappa$ and $N\geq2$. Write $x_0=c(R)$, $\ell=\ell(R)$, and let
\[
f_N\coloneqq\frac{(\cdot-x_0)^\gamma}{\gamma!}\Omega_N,
\qquad \Omega_N\coloneqq\sum_{j=1}^N\omega^\infty_{Q_j},
\qquad c(Q_j)=x_0,\quad\ell(Q_j)=2^j\ell.
\]
These are the dilated tests of \cite{DPModPoly}*{\S11.3}. On $\frac32Q_1\supseteq R$,
$\partial^{\gamma'}f_N=\delta_{\gamma'\gamma}N$ for $|\gamma'|=\kappa$, so the last sum in \eqref{e:SoTcoef} equals $Nb^{\alpha,\gamma}_R/\gamma!$.
For each $|\beta|=\kappa$, expanding the derivatives of the cutoffs gives
\[
\left|\partial^\beta f_N(x)-\delta_{\beta\gamma}N\right|
\lesssim1+\log\left(1+\frac{|x-x_0|}{\ell}\right),
\]
uniformly in $N$. Each term of this expansion is bounded and vanishes on $\frac32Q_j$, so at most $O(1+\log(1+|x-x_0|/\ell))$ terms contribute. The vanishing mean and decay of $\varphi_R^{\alpha,\beta}$ therefore give
\[
\left|\varphi_R^{\alpha,\beta}(\partial^\beta f_N)\right|\lesssim1,
\qquad
\mathrm S^\tau(\partial^\alpha Tf_N)
\gtrsim\left(\frac N{\gamma!}|b_R^{\alpha,\gamma}|-C\right)_+\ind_R.
\]
Use the normalized power weight
\[
w_N(x)\coloneqq\ell^{-d}\left|\frac{x-x_0}{\ell}\right|^{-d+1/N}.
\]
By \cite{DPModPoly}*{Lemma 11.5}, $[w_N]_{A_\infty}\sim N$ and $w_N(R)\sim N$. Also, Lemma 11.4(a) there, as $\mathcal N_1$ commutes with dilations, gives
\[
\mathcal N_1(\partial^\beta f_N)(x)
\lesssim\frac1r\min\left\{1,\left(\frac{2^N\ell}{|x-x_0|}\right)^d\right\},
\qquad
\|\mathcal N_1(\partial^\beta f_N)\|_{L^2(w_N)}\lesssim\frac{\sqrt N}{r},
\]
the norm estimate following by polar integration. Applying \eqref{e:SoT1} at $p=2$,
\[
\left(\frac N{\gamma!}|b_R^{\alpha,\gamma}|-C\right)_+\sqrt N
\lesssim\frac{C_2}{r}N.
\]
Divide by $N^{3/2}$ and let $N\to\infty$. Thus $b_R^{\alpha,\gamma}=0$.
\end{proof}

\subsection*{Declaration of generative AI and AI-assisted technologies
in the manuscript preparation process}
While this work was being finalized,  the author used Claude
(Anthropic) and ChatGPT (OpenAI) to assist in
checking the mathematical arguments, cross-references and citations,
and in preparing and editing the \LaTeX{} source. The author reviewed
and edited the manuscript, independently verified the mathematical
arguments and references, and takes full responsibility for the
content of the article.

\bibliographystyle{plain}
\bibliography{cancellativesparse}
\end{document}